\documentclass[12pt]{article}
\usepackage{amsmath,amssymb,amsbsy,amsfonts,amsthm,
latexsym,amsopn,amstext,amsxtra,euscript,amscd,amsthm}

\usepackage{relsize}
\usepackage{enumitem}
\usepackage{graphicx}
\usepackage{float}
\usepackage{wrapfig}
\usepackage{picinpar}
\usepackage{tcolorbox}
\usepackage{hyperref}
\newcommand{\p }{p_n\# }
\def\bsq{\blacksquare}
\usepackage{tikz}

 \newif\ifstartedinmathmode
\newcommand\encircled[1]{%
  \relax\ifmmode\startedinmathmodetrue\else\startedinmathmodefalse\fi%
  \tikz[baseline,anchor=base]{%
  \node[draw,circle,outer sep=0pt,inner sep=.2ex]
    {\ifstartedinmathmode$#1$\else#1\fi};}%
}

\title{\bf A special sequence and primorial numbers }
\author{
{\sc Amit Kumar Basistha}\\
{\rm Indian Statistical Institute}\\
{ Bangalore Centre }\\
{8th Mile, Mysore Rd, RVCE Post}\\
{ Bengaluru, Karnataka 560059}\\
{\tt basisthaamitkumar2@gmail.com} \\
{\sc Eugen J. Ionascu}\\
{\rm Columbus State University}\\
{4225 University Avenue}\\
{Columbus, GA 31907}\\
{\tt ionascu\_eugen@columbusstate.edu} \\
 }
\begin{document}
%MSC: 11A41;11A51;11B39;11N05
\date{ }
\maketitle

%\baselineskip=1.7\baselineskip

\begin{abstract}
In this paper, we study a class of functions defined recursively on the set of natural numbers in terms of the greatest common divisor algorithm of two numbers and requiring a minimality condition.  These functions are permutations, products of infinitely many cycles that depend on certain breaks in the natural numbers involving the primes, and some special products of primes with a density of approximately $29.4\%$. We show that these functions split into only two equivalence classes (modulo the natural equivalence relation of eventually identical maps): one is the class of the identity map and the other is generated by a map whose discrete derivative is almost periodic with ``periods"  the primorial numbers.
\end{abstract}

%%%%%%%%%%%%%%%%%%%%
%%  PAPER BEGINS  %%
%%%%%%%%%%%%%%%%%%%%

\def\RR{{\rm I}\!{\rm R}}
\def\fp#1{(#1)}

\newtheorem{thm}{Theorem}
\newtheorem{rem}{Remark}
\newtheorem{defn}{Definition}
\newtheorem{lem}{Lemma}[section]
\newtheorem{theorem}[lem]{Theorem}
\newtheorem{alg}[lem]{Algorithm}
\newtheorem{cor}[lem]{Corollary}
\newtheorem{conj}[lem]{Conjecture}
\newtheorem{prop}[lem]{Proposition}
\newtheorem{heu}[lem]{Heuristic}
\newtheorem{LMN}[lem]{Theorem}
\newtheorem{exa}{Example}

\def\blacksquare{{\ \vrule height7pt width7pt depth0pt}}
\def\bsq{\blacksquare}
\def\go{\rightarrow}
\def\ds{\displaystyle}
\def\n{\noindent}

%%%%%%%%%%%%%%%%%%%%%%%
\let\workingver=n
\def\begeq#1{\begin{equation}\mylabel{#1}}
\def\endeq{\end{equation}}
\def\begalg{\begin{alg}}
\def\endalg{\end{alg}}
\def\be{\beta}
\def\be{\beta}
\def\de{\delta}
\def\st{\star}
\def\op{\oplus}
\def\Llr{\Longleftrightarrow}
\def\Om{\Omega_f}
\def\rs{{\it RotS} }
\def\hF{{\hat {\cal F}}}
\def\refeq#1{\if\workingver y(\ref{#1})-
[[#1]]\else(\ref{#1})\fi}
\def\refth#1{\if\workingver y\ref{#1}-[[#1]]\else\ref{#1}\fi}
\def\mylabel#1{\if\workingver y\label{#1}{\bf\ \ [[#1]]\ \ }
\else\label{#1}\fi}
\def\mybibitem#1{\if\workingver y\bibitem{#1}{\bf\ \ [[#1]]\ \
}
\else\bibitem{#1}\fi}
%%%%%%%%%%%%%%%%%%%%%%%%%%%%%%%%%%%%%%%%%%%
\newcommand{\kro}[2]{\left( \frac{#1}{#2} \right) }
%%%%%%%%%%%%%%%%%%%%
%%  PAPER BEGINS  %%
%%%%%%%%%%%%%%%%%%%%

\section{Introduction} \label{intro}

\n  The following problem was proposed by the first author in Crux Mathematicorum\cite{amit} :

  \textit{ Let $f:\mathbb N\to \mathbb N$ with $f(1) = 1$,   $f(2) = a$ for some $a\in \mathbb N$ and, for each positive integer $n\ge 3$, $f(n)$ is the smallest value not assumed at lower integers that is
coprime with $f(n-1)$. Prove that $f$ is onto.  }

\n In what follows we are going to use the conventional notation for the greatest common divisor of two natural numbers: for $m,n\in \mathbb N$ this will be denoted by $\gcd(m,n)$.

If $a$ is not relevant we will simply refer to the sequence by $f$, but sometimes $a$ may play an important role in which case we will use $f_a$ instead. For example, $f_1$ is not much different from $f_2$ which is the identity map on $\mathbb N$, so we will assume that $a\ge 2$ from now on.

A special case of the sequence was introduced for $a=3$, in OEIS (The On-Line Encyclopedia of Integer Sequences) by Reinhard Zumkeller in 2003 as \href{https://oeis.org/search?q=+A085229&sort=&language=&go=Search}{A085229} but with an intrinsic definition: ``Smallest natural number $x_n$ which is coprime to $n$ and to $x_{n-1}$, and is not yet in the sequence ($x_1=1$).'' It is not obvious that this definition is equivalent to our definition above and we will show this as a corollary of Theorem~\ref{theorem1}. A little different definition is given for A123882 which coincides with our sequence $f_3$, for all indices $n\ge 4$.

It turns out that the sequence can be easily computed (a short code in Python is provided on OEIS) and the first \href{https://oeis.org/A085229/b085229.txt}{10,000 terms} of $f_3$ are also available. We will prove that $f_a$ is not only a surjection but also an injection. So, we are dealing with permutations of
$\mathbb N$. In particular, A085230 is $f_3^{-1}$.

There are a few results that are mentioned in OEIS by Michael De Vlieger (April 13th, 2022) concerning properties of $f_3$ (see \cite{mdv}).

\begin{prop}
\label{prop:1}
  $f_3(2k+1)=2k$ for all $k>0$. 
\end{prop}

We will show that this follows from Theorem~\ref{theorem1}.

\begin{prop}
\label{prop:2}  $f_3(3k+1)=3k$ for all $k>1$. 
\end{prop}

Let us observe that for $k=2m$ this follows from Proposition~\ref{prop:1}. Also, from Proposition~\ref{prop:1}  we see that the terms of $f_3$ are following the pattern:

$$1,3,2,5,4,7,6, \  \boxed{ ? }\ , 8,  \ \boxed{? } ,\   10, \ \boxed{? }  12, \  \boxed{ ?} ,\  14, \boxed{ ?}  16,  \ \boxed{ ?} , \ldots$$

Let us assume for the moment that $f_3$ is a surjection. Then, if we look at $9$, it cannot fit in the first box since $\gcd(6,9)=3$, and so it should go into the second by its minimality. That means $f(10)=9$. Then $15$ cannot go into the first, the third, or the fourth box since
$\gcd(12,15)=3$ and then by minimality it has to go into the fifth which means $f(16)=15$. This argument can be finished by induction showing that $f(6m+4)=6m+3$ for $m\ge 1$ proving the claim for $k=2m+1$.

The list of the first $24$ terms in $f_3$ is included next:
 $$\begin{array}{ccccccccccccc}
 n& 1 & 2 & 3 & \boxed{4} & 5 & \boxed{6 }& 7 & \boxed{8} & 9 & 10 &11& \boxed{12} \\
  \hline
f(n)&  1 & 3 & 2 & \encircled{5} & 4 &  \encircled{7} & 6 &  \encircled{11} & 8 & 9 &10& \encircled{13}
\end{array}$$

$$\begin{array}{ccccccccccccc}
 n& 13 & \boxed{14} & 15 & \boxed{16} & 17 & \boxed{18 }& 19 & \boxed{20} & 21 & 22 &23 & \boxed{24} \\
  \hline
f(n)&  12 &  \encircled{17} & 14 & 15 & 16 &  \encircled{19} & 18 &  \encircled{23} & 20 & 21 & 22&  \encircled{25}
\end{array}$$

In \cite{mdv} there is a mention of a concept named {\bf record}. If we look at the above table, we observe some bigger jumps when the sequence goes up more than $1$ from the previous value in the sequence. In the next section, we will introduce a slightly different term, that of a {\bf turning point}. The numbers in the boxes are turning points and their values (circled) are records. The smallest composite value for a turning point is $f_3(24) = 25$ and the smallest record that has at least two prime factors is $f_3(54) = 55$. From Proposition~\ref{prop:1} and Proposition~\ref{prop:2}, we see that every record is an odd number and $3$ cannot divide a record. Hence, all records must be of the form $6k\pm 1$. We will show that every prime $p\ge 5$ is a record. So, this function contains good information about primes having the advantage that
the terms can be recursively calculated only using the gcd function. Also, one can compute a section of the sequence without knowing all of the terms up to that particular starting input.

Here is a list of all the non-prime records less than $100$ and their jumps, i.e., $j_r=r-f_3^{-1}(r)$.

$$ \{ [25, 1],\  [49, 1],\ [55, 1], \ [77, 3],\ [85,1],\ [91,1]\}$$

\n There are a few important questions here related to the records (especially the ones which are composite numbers), say $\{\overline{R}_j\}$, $\overline{R}_1=25$, $\overline{R}_2=49$, $\overline{R}_3=55$, etc. What is their distribution? What is their distribution within the records, or equivalently, what is the distribution of the primes within the set of records?

In general, a permutation of a finite set is a product of cycles. In our case, $f_a$ is a permutation of the infinite set $\mathbb N$. However,  we will show that $f_a$ is still a product of finite cycles. We use the usual convention of denoting a cycle by $(c_1,c_2,\ldots,c_n)$ meaning the permutation which maps $c_1$ into $c_2$, $c_2$ into $c_3$,$\ldots$, and $c_n$ into $c_1$. Cycles of length one are usually left out.
This way we can write
$$f_3=(3,2)(5,4)(7,6)(11,10,9,8)(13,12)(17,16,15,14)(19,18)(23,22,21,20)(25,24) \ldots$$
So, essentially $f_3$ is defined by the sequence of records.

\begin{figure}[H]
    \centering
    \includegraphics[height=0.4\textwidth,width=0.6\textwidth]{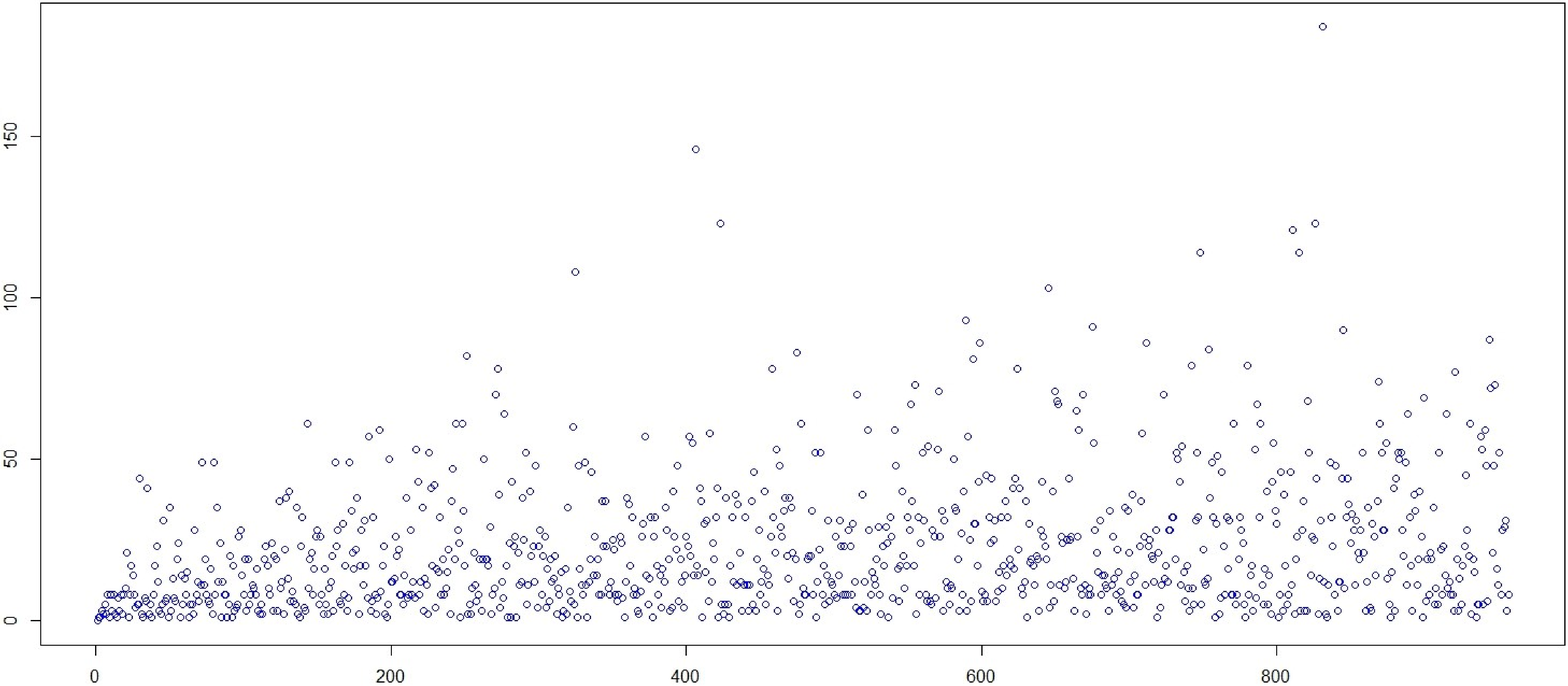}
    \caption{Twin primes distribution into cycles}
\end{figure}

We can permute these cycles in whatever order we like. But we will assume the set of these cycles to be ordered in terms of the numbers in it (non-decreasing). For each record $t$,  denote by $C(t)$ the cycle number. For instance, $C(23)=8$ and $C(25)=9$.

In Figure~1, we included the values  of $C(m_{j+1})-C(M_j)$ where $(m_j,M_j)$ is a twin pair of primes.
We will formulate a few conjectures about this data in Section 3. We notice that $j(r)=1$ is an indication that $r$ is the biggest of a twin pair. This happens for a lot of records which are not primes. However, let us call these records {\bf twin records}. 

Let us point out at least one connection with primorial numbers. If $p_n$ is the $n$-th prime, then the $n$-th primorial number (see \cite{pm})   is defined by
$$p_n\# =\prod _{k=1}^n p_k. $$

\n The values of $p_n\# $ for $n=1$, 2, $\ldots$, are 2, 6, 30, 210, 2310, 30030, 510510, $\ldots$ (\href{https://oeis.org/search?q=A002110&sort=&language=&go=Search}{OEIS A002110}).

In Figure~2, we included the values of $g(t):=f(t)-f(t-1)$, $t=1,2,\ldots,12000$ for $f_3$.

\begin{figure}
    \centering
    \includegraphics[height=0.4\textwidth,width=0.6\textwidth]{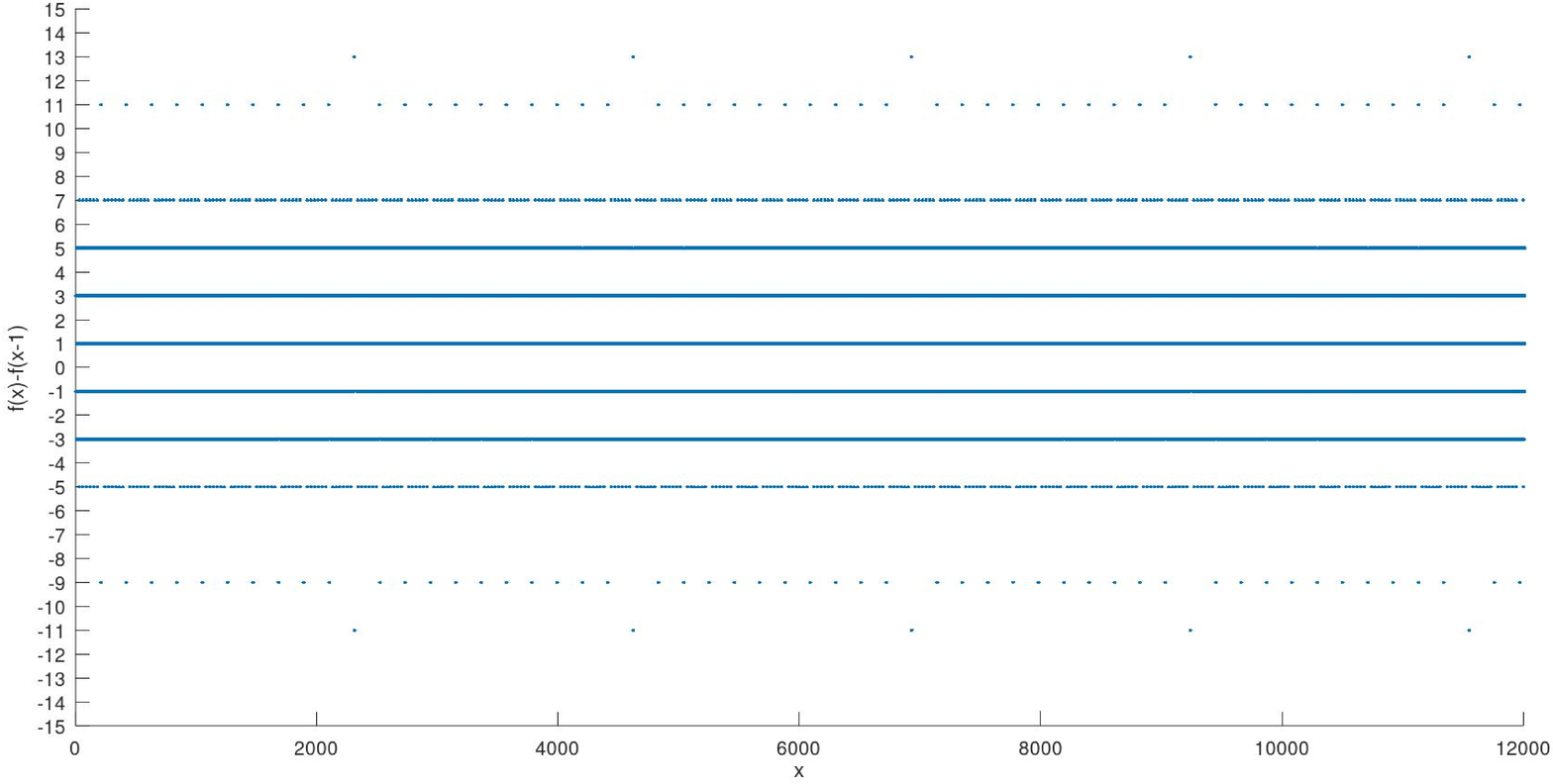}
    \caption{Discrete Derivative of f}
\end{figure}

We will show that $|g(t)|$ is unbounded by proving that for infinitely many $k\in \mathbb N$, we have
$$g(k p_n\#+1)\ge 2n+1.$$
\par

A similar sequence was studied in \cite{Apple} (2015) This is sequence \href{https://oeis.org/A098550}{A098550} in the OEIS. Very similar results are shown, including the proof of the one-to-one correspondence with the natural numbers.

\section{Results and proofs}

\n \begin{defn} A {\color{blue} \it  turning point} is a natural number having either of the following two properties:
\begin{enumerate}[label=\roman*)]
    \item $t>3$ and $f(t)-f(t-1)>1$.
    \item $t=3$ and $f(3)\not = \min \mathbb N\setminus \{1,a\}$.
\end{enumerate}
The value $f(t)$ for $t$ a turning point is called a {\color{blue} \it record}.\end{defn}
For example if $a=4$, the sequence goes like
$$1,4,3,2,5,6,7,8, \ldots $$
and so $3$ is a turning point and $t=5$ is also a turning point since $f_4(5)-f_4(4)=5-2=3>1$. Since $f_4(5)=5$, this is called a fixed point and we see that lots of fixed points follow, making the sequence less interesting. In fact, for $a=2$, we have no turning point since $f_2(n)=n$ for all $n$.
We are not going to consider these functions and assume further that $a\ge 3$.

\n \begin{defn} An {\color{blue} \it essential turning point}  (or ETP) is a turning point $t$ having in addition the following three properties: 
\begin{enumerate}[label=\roman*)]
    \item  $t>a$ and $f(t)\not =t$.
    \item   $f(t-1)=t-2$.
    \item $\{1,2,3,\ldots,t-1\}=\{f(1),f(2),\ldots,f(t-1)\}.$
\end{enumerate}
\end{defn}

For example, if $a=7$, the first 12 terms are given in the table below:

$$\begin{array}{ccccccccccccc}
 n& 1 & 2 & 3 & 4 & 5 & 6 & 7 & \boxed{8} & 9 & 10 &11& \boxed{12} \\
  \hline
f(n)&  1 & 7 & 2 & 3 & 4 & 5 & 6 & \encircled{11} & 8 & 9 &10&\encircled{13}
\end{array}$$

We observe that $t_1=8$ is an essential turning point. Also, $t_2=12$ is an ETP  and the list continues. The two corresponding records are the primes $11$ and $13$. We prove the following theorem about ETPs.

%%%%%%%%%%%%%%%%%%%%%%%%%%%%%%%%%%%%%%%%%%%%%%%%

\begin{theorem}\label{theorem1}   If $t$ is an ETP then $T:=f(t)+1$ is the next  ETP and
 there are no turning points in the interval $(t,T)$.
\end{theorem}
%%%%%%%%%%%%%%%%%%%%%%%%%%%%%%%%%%%%%%%%%%%%%%%%

\n \begin{proof}  Since $t>a$,  we may assume that $t>3$, and define $p:=f(t)-f(t-1)$. Hence, we have $p>1$ since $t$ is a turning point. Hence, $T=f(t)+1=p+f(t-1)+1=p+t-2+1=p-1+t>t$.

By the definition of $f$, we have $1=(f(t),f(t-1))=(f(t-1)+p,f(t-1))=(p,f(t-1))$. So $p$ is the smallest natural number $p>1$, with the property $(p,f(t-1))=1$, such that $p+f(t-1)=p+t-2\ge t$ is not one of the values  $f(m)$ with $m\le t-2$.  The last condition is obviously satisfied since $\{f(1),f(2),\ldots,f(t-2)\}\subset \{1,2,\ldots,t-1\}$ by the definition of ETP.

We claim that $(f(t),f(t-1)+2)=1$. By way of contradiction, suppose that  $$(f(t),f(t-1)+2)=d>1.$$ Then $d$ divides $f(t)-(f(t-1)+2)=p-2<p$. But $(f(t),f(t-1))=1$ and since $d$ divides $f(t)$ we must have $1=(d,f(t-1))=(d+f(t-1),1)$. The minimality of $p$ that was pointed out above shows that $d=p$, a contradiction. It remains that $d=1$ and then by the definition of $f(t+1)$,  $f(t+1)=f(t-1)+2=t-2+2=t$.

As result of this, $f(t+2)=t+1$, \ldots, $f(t+j)=t+j-1$, as long as $t+j-1<f(t)=T-1$.
None of the values $t+j$ are turning points. For $j$ such that  $t+j-1=f(t)=T-1$ or $j=j_0=T-t$,  we have

$$\{f(1),f(2),\ldots,f(t-1),f(t),\ldots,f(T-1)\}=\{1,2,\ldots,t-1,T-1,t,t+1,\ldots,T-3,f(T-1)\}=$$
$$\{1,2,\ldots,t-1,t,t+1,\ldots,T-2,T-1\}$$
which shows that $T$ is the next  ETP provided that  $f(T)-f(T-1)=f(T)-T+2>1$. Because $f(T)$ is forced to be more than or equal to  $T$ the last constraint is satisfied.  Therefore, the next ETP is $T$.
\end{proof}

\n \begin{rem} If we have at least one  ETP, say $t_0$ (may as well assume it is the smallest one),  then we can generate them all by using the recursion 
$$t_n=f(t_{n-1})+1\ for \ n\ge 1.$$ \end{rem}

The sequence $t_n$ is strictly increasing and so it is unbounded. The property (iii) of an ETP shows that $f_a$ is then onto. So, the problem we started with in the Introduction is proven if we show the existence of at least one ETP. In general, for some values of $a$, $f_a$, doesn't have any $ETP$. However, in that case, it will be easy to show that $f_a$ is a bijection.

\n \begin{rem} Let us observe that $f$ is actually  one-to-one. Indeed, let us assume that $1\le m<n$. If $n=2$ then $m=1$ and so $f(1)=1\not=f(2)=a$ by the assumption in the Introduction. If $n\ge 3$,  by definition $f(n)$ is not in the set
$$\{f(1),f(2),\ldots,f(m),\ldots,f(n-1)\}$$

\n and so $f(n)\not=f(m)$.  We have seen that
$$f_3=(3,2)(5,4)(7,6)(11,10,9,8)(13,12)(17,16,15,14)(19,18)(23,22,21,20)(25,24) \ldots$$
\n with the beginning of each cycle an ETP (except for $2$). Similarly, we have the writing
$$f_7=(7,6,5,4,3,2)(11,8,9,10)(13,12)\ldots $$
\n and clearly we observe that $f_3$ and $f_7$ have the same cycles eventually.
Also, another interesting situation appears if $a$ is a multiple of $6$:

$$f_6=(6, 5, 2, 3, 4)$$
$$f_{12}=(12, 5, 2, 3, 4, 7, 6, 11, 8, 9, 10)$$
$$f_{18}=(18, 5, 2, 3, 4, 7, 6, 11, 8, 9, 10, 13, 12, 17, 14, 15, 16),\ldots $$
in which case $f(n)=n$ eventually (for big enough $n$), and $f$ consists of only one nontrivial cycle.
\end{rem}

\n \begin{defn}\label{Def 3}  Let's call two permutations $f_a$ and $f_b$ {\color{blue} \it  EI-permutations}  (eventually identical) if there exists $m$ which depends on $a $ and $b$ such that $f_a(n)=f_b(n)$ for all $n>m$. This (equivalence) relation partitions the set of these bijections into equivalence classes, $\cal C$.\end{defn}
It looks like we have only two classes so let $\color{red} {\cal C}={\cal C}_3 \cup {\cal ID}$ where ${\cal C}_3$ is the class of $f_3$ and ${\cal ID}$ is the class of $f_2$ or eventually the identity maps. For $a$ such that $f_a\in {\cal ID}$ let us denote by $M_a$ the smallest natural number with the property $f_a(n)=n$ for all $n\ge M_a$. Also, we will use the notation
$${\cal A}:=\{a\in \mathbb N| f_a \ \ \ \ \text{is in}  \ \ \ {\cal ID}\}.$$

The set ${\cal A}$ appears to be nontrivial but we will see that most of the numbers which are multiples of $6$ are in ${\cal A}$. There are some exceptions such as $a=216=210+6$. We observe that $210=2\cdot 3\cdot 5\cdot 7$ and $6=2\cdot 3$ which are primorial numbers.
We prove next Proposition~\ref{prop:1} from the Introduction and that every prime is a record for $f_3$.

\begin{cor}\label{corollary1}   $f_3(2k+1)=2k$ for all $k\geq 1$ and if $p$ is a prime greater than or equal to $5$, $p$ is a record of $f_3$.
\end{cor}

\n \begin{proof} It is clear that the first ETP for $f_3$ is $t_1=4$, the second is $t_2=f(t_1)+1=6$, and so on.  We see that between two consecutive turning points, as in the above proof, the sequence continues as :

$$\begin{array}{cccccccc}
 k  & t_n-1  &t_n      & t_n+1 & t_n+2, \ldots &  T_n-1      &  T_n         & T_n+1 \\
  \hline
f_3(k)&  t_n-2 &T_n-1 & t_n    & t_n+1,\ldots,& T_n-2& f(T_n) & T_n,\ldots
\end{array}, \ \ T_n=f(t_n)+1.$$
Using induction on $n$ we see that every $t_n$ must be even and every record $R_n=f(t_n)$ must be odd. So the sequence of values that are even goes in increasing order and $f_3(2k+1)=2k$ for each $k\ge 1$.

\n The second claim in the corollary is obviously true for $p=5$. Let us assume by way of contradiction that $p> 5$ is not a record. So, it will appear in the sequence ($f_3$ is a bijection) between two EPTs as above. But then $p<T_n-1=f(t_n)$ and $p$ is relatively prime with $t_n-2<p$, contradicting the choice of $f(t_n)$.Then every prime (except $2$ and $3$) is a record.
\end{proof}

 \begin{rem} We see that the proof above works if we assume that for the prime $p$ there exists an ETP $t$ such that $p>t$.  In particular, it is true for every $a$ such that $f_a$ has at least one ETP and $p$ a prime big enough.
\end{rem}

Next, let us show that $|g(t)|$ is unbounded, where $g(t)=f(t+1)-f(t)$ for $t\in \mathbb N$.

\begin{prop}\label{prop:3}$g(k p_n\#+1)\ge 2n+1$ for infinitely many $k\in \mathbb N$.
\end{prop}

\begin{proof} By Dirichlet's theorem on arithmetic progressions, $q=k p_n\#+1$ is a prime for infinitely many $k\in \mathbb N$. Choose $k$ so that $q>5$ . This means that $q$ is an ETP for $f_3$ by Corollary~\ref{corollary1}. Then we have
$$g(k p_n\#+1)=f_3(k p_n\#+2)-f_3(k p_n\#+1)=f_3(k p_n\#+2)-k p_n\#:=m+1-q$$
where $m:=f_3(k p_n\#+2)$ is relatively prime with $q-1=k p_n\#$ and is bigger than $q$. Then $m$ must not be divisible by any of the prime factors of $q-1$. So $m\neq q+1$ as $q+1$ is even. So we are done if $n=1$.

Suppose, $n\ge 2$. Observe that $\gcd(q-1,q+2r)$ is divisible by at least one of the $p_i$, $i=1,2,\ldots,n$ due to the obvious inequality $p_n\ge 2n-1$. Also as $q-1$ is even $m$ must be odd and hence $m-q$ is even. So $m-q\geq 2n$ which completes the proof.
\end{proof}

\begin{theorem}\label{theorem2}   If $a$ is odd and $t$ is an ETP then  $t$ is an even number and $f_a\in {\cal C}_3$
\end{theorem}

\begin{proof}  If $a$ is an odd number then $f(3)=2$, $f(4)=3$, and so on, until $f(a)=a-1$, and then $f(a+1)$ is $a+2$ or bigger,  turning $a+1$ into an EPT. In fact, this is the first EPT. Clearly, $t_1=a+1$ is even and  so $f(t_1)$ must be odd, otherwise $(f(a),f(a+1))=(a-1,f(t_1))\ge 2$, a contradiction.  This shows that $t_2=f(t_1)+1$ (by Theorem~\ref{theorem1}) is even. Inductively, we see that all of ETP's must be even.
Supposed that we take a prime $p>t_2=f(t_1)+1$ and also $p\ge 5$. As in the proof of Corollary~\ref{corollary1} $p$ must be a record or  $f(t_k)=p$ for some $k$. Then $t_{k+1}=f(t_k)+1=p+1$ is a EPT for $f_a$ but also for $f_3$. Therefore from this point on $f_3(n)=f_a(n)$ for all $n\ge p+1$ since the definitions of the two functions are recursively in terms of the same data.
\end{proof}

\par
\begin{rem} The result obtained in Theorem~\ref{theorem2} can be clearly improved by only assuming that $f_a$ is a function which does have an ETP.
\end{rem}

%%%%%%%%%%%%%%%%%%%%%%%%%%%%%%%%%%%%%%%
\begin{theorem}\label{theorem3}  Assume {\bf $ a$} is even and a multiple of $6$. Then $f_a$ is either in $\cal ID$ or in ${\cal C}_3$ in which case every ETP is even.
\end{theorem}
%%%%%%%%%%%%%%%%%%%%%%%%%%%%%%%%%%%%%%%%%%

\begin{proof}  Let us assume that $a>4$ and define $\kappa:=f(3) \not =3 $ which must be an odd number and in addition $(\kappa,a)=1$.  Since $a-1$ is odd and $(a-1,a)=1$, by the minimality of $f(3)$ we see that $\kappa\le a-1$.

Clearly then $\kappa$ is an odd number greater than or equal to $5$  and then $f(4)=2$, $f(5)=3$, \ldots, $f(\kappa+1)=\kappa-1$ an even number. Then, $f(\kappa+2)$  should skip $\kappa$ since it is already in the list. Then the next candidate is $\kappa+1$ but this is also even so we need to move up to $\kappa+2$, i.e., $f(\kappa+2)=\kappa+2$. This is possible for lots of values of $a$, ($a\in \{6,12,18,24,36,\ldots\}$). If we have
$$\{1,2,\ldots,\kappa+2\}=\{f(1),f(2),\ldots,f(\kappa+2)\}$$
then clearly $f(n)=n$ for all $n\ge \kappa+2$. In this case, we have no turning point and no ETP.
This is exactly what happens if $a=6$ and only if $a=6$ (but it is not necessary to prove this at this point). Hence, we will assume that $a\ge 12$ from here on.  So, if $f(\kappa+2)=\beta$ for some odd number $\beta\ge \kappa+2$ with $(\beta,\kappa-1)=1$, $\kappa+2$ becomes a turning point. We will look at $f_{36}$ as a generic example:

$$\begin{array}{ccccccccccccccccccc}
 n& 1 & 2 & \boxed{3} & 4 & 5 & 6 & \boxed{7} & 8 & \boxed{9} & 10 &11& 12&\boxed{13}& 14& \boxed{15}& 16 &17& 18 \\
  \hline
f(n)&  1 & 36 & 5 & 2 & 3 & 4 & 7 & 6 & 11 & 8 & 9&10& 13& 12&  17&14 &  15& 16
\end{array}$$

$$\begin{array}{ccccccccccccccccccc}
 n& \boxed{19} & 20 & \boxed{21} & 22 & 23 & 24 & \boxed{25} & 26 & \boxed{27} & 28 &29& 30&\boxed{31}& 32& \boxed{33}& 34 &35& 36 \\
  \hline
f(n)&  19  & 18 & 23 & 20 & 21 & 22 & 25 & 24 & 29 & 26 & 27 &28 & 31 & 30&  37& 32 &  33& 34
\end{array}$$
$$f(37)=35\ \ \text{and}\ \ \ f(38)=38$$

This shows that $M_{36}=38$. We notice that in this example $f(2k)=2k-2$ for all $k\in \{2,3,\ldots,18\}$.
This rule breaks at $k=19$ and also $38$ is a turning point and the last one.

Let us denote by $s$, the largest turning point with the property that $s<a$ and we let $q=f(s)$.

 We claim that if $q+1=a$ then $q+2=a+1$ is a turning point which is equal to $M_a$. If $q+1>a$, then $q+1$ is a turning point  which  is  equal to $M_a$ or it is the first ETP  and $f_a$ is in ${\cal C}_3$

In general for $a\in \{6,12,\ldots\}$, we already know that $f(4)=2$, $f(5)=3$, $f(6)=4$, and so on.
First, let us prove that

\begin{equation}\label{eq1} f(2k)=2k-2\ \ \text{ for all }\ \  k\ \ \text{such that} \ 4\le 2k< q+1
\end{equation}
   We have already shown that ($a$ is at least $12$)
$f(\kappa+1)=\kappa-1$ and $f(\kappa+2)=\beta\ge \kappa+2$.
We claim that $f(\kappa+3)=\kappa+1$ which means the property (\ref{eq1}) holds for all $k$ such that $2k\le \kappa+3$. The list
$$[f(1),f(2),\ldots,f(\kappa+2)]=[1,a,\kappa,2,3,4,\ldots,\kappa-1,\beta]$$
shows that if $$(\kappa+1,\beta)=1$$ then $f(\kappa+3)=\kappa+1$. The proof of this is done by way of contradiction as in the proof of Theorem~\ref{theorem1}. Hence, the property (\ref{eq1}) holds for $2k=\kappa+3$. This allows us to continue the above list
$$[f(1),f(2),\ldots,f(\kappa+2),\ldots,f(\beta+1)]=[1,a,\kappa,2,3,4,\ldots,\kappa-1,\beta,\kappa+1,\ldots,\beta-1]$$

\n making $\beta+2$ the next turning point. This list continues until we get to the last turning point less than $a$, which is $s$.  Note that all the turning points are odd numbers. The list above becomes

$$[f(1),f(2),\ldots,f(\kappa+2),\ldots,f(\beta+1),f(\beta+2),\ldots,f(s),\ldots,f(q+1)]=$$
$$=[1,a,\kappa,2,3,4,\ldots,\kappa-1,\beta,\kappa+1,\ldots,\beta-1,\beta',\ldots,q,\ldots,q-1]$$

We observe that $q+2<a$ is not possible since $s$ is the greatest turning point less than $a$. So, $q+2\ge a$, but since $q+2$ is odd we must have $q+2\ge a+1$. If we have equality, $q=a-1$, then from the above equality of lists we conclude that $f(q+2)=a+1=q+2$ ($\gcd(q+2,q-1)=\gcd(3,q-1)=\gcd(3,a-2)=1$)   and for $n\ge q+2$, we have $f(n)=n$. That makes $a+1=M_a$ and there are no EPT.

In the second situation, $q+1>a$ and the above equality of lists need to be corrected. So

$$[f(1),f(2),\ldots,f(\kappa+2),\ldots,f(\beta+1),f(\beta+2),\ldots,f(s),\ldots,f(q),f(q+1)]=$$
$$=[1,a,\kappa,2,3,4,\ldots,\kappa-1,\beta,\kappa+1,\ldots,\beta-1,\beta',\ldots,q,\ldots,q-2, q+1]$$
\n because in the process of writing the even numbers we had to skip over $a$ if  $\gcd(q-2,q+1)=1$.
This will place $f_a$ in $\cal ID$. If $\gcd(q-2,q+1)>1$, we get an ETP and we have $f_a\in {\cal C}_3$.
\end{proof}
%%%%%%%%%%%%%%%%

\begin{rem} In the last proof if $a=216$ then in the last step  $q=221$ and $\gcd(q-2,q+1)=3$. That places $f_{216} \in {\cal C}_3$.
\end{rem}

%%%%%%%%%%%%%%%

Let us denote by $R$ the set of records for $f_3$.

\begin{theorem}\label{theorem4}  The set    ${\cal A}$ has the following description

\begin{equation}\label{descriptionofA}
\mathcal{A}=\{2,4\}\cup \{a| a\equiv 0 \pmod{6} \ \text{such that there exists a record}\  r\in R, \  \ |r-a|\le 1\}
\end{equation}

\end{theorem}

\begin{proof} We are denoting the right side of (\ref{descriptionofA}) by $B$. First we show that $B\subset {\cal A}$.
We have already observed that $2$, $4$ and $6$ are in  ${\cal A}$.  Consider $p$ the smallest prime   less than $a$
such that $p$ doesn't divide $a$ ($\gcd(p,a)=1$). This prime exists since $a\ge 12$.  Then $f_a(3)\le p$  and $3$ is a turning point for $f_a$. Using the same arguments as before, we see that eventually the records of $f_a$ are going to overlap with the records of $f_3$ ($p$ in particular has to be a record of $f_a$ also). Because $p$ is less than $a$, we may assume that this happens before $a$, i.e., the records of $f_a$ which are also records of $f_3$ start at a value less than $a$ (this doesn't happen if $a=6$).  Then it makes sense to define $L$ be the largest record (in $R$) less than $a$ and $S$ be the smallest one (in $R$) bigger than $a$.  We observe that $L$ and $R$ are records for $f_a$ too (by earlier observations). Let us analyze the two cases.

We assume first that $L=a-1$. The sequence $f_a$ takes the values
$$1,a,\ldots,L,k,k+1,k+2,\ldots,L-1=a-2 \ \ \ \text{where} \ f(k-1)=k$$
\n and by the maximality of $L$, we must have used all of the values in the set $\{1,2,3,\ldots, a\}$. Therefore, we must have  $f(a)=a-2$ and then $f(a+1)=a+1$  ($\gcd(a+1,a-2)$ can be at most $3$ and since $6$ divides $a$, $3$ doesn't divide $a-2$). Now it is easy to see that $f_a\in {\cal ID}$.

We assume next that $L<a-1$, and $S=a+1$ ($a=36$ is the first with such a property).  Then the
sequence $f_a$ takes the values
$$1,a,\ldots,L,k,k+1,k+2,\ldots,L-1<a-2,S,m,m+1,\ldots,S-2=a-1 $$

 \n where $f(k-1)=k$  and $f(m-1)=S$, which cover all the values in the set $\{1,2,\ldots,a,S\}$. Therefore, we must have $f(a+1)=a-1$ and so $f(a+2)=a+2$ ($\gcd(a+2,a-1)$ can be at most $3$ and since $6$ divides $a$, $3$ doesn't divide $a-1$). Now it is easy to see that $f_a\in {\cal ID}$. This shows that $B\subset {\cal A}$.

Suppose that $a\not \in  B$. First we assume that $a$ is such that
$L<a-1$ and $S>a+1$. Because $a$ is even, $\gcd(a-1,a+1)=1$ and so   the
sequence $f_a$ takes the values
$$1,a,\ldots,L,k,k+1,k+2,\ldots,L-1<a-2,S,m,m+1,\ldots,a-1,a+1,\ldots, S-1$$
which cover all the values in the set $\{1,2,\ldots,a,\ldots, S-1,S\}$. This makes $f(S)=S-1$ and so $S+1$ is an  ETP placing  $f_a\in {\cal  C}_3$ and so $a\not\in {\cal A}$.

Next, we assume that $a$ is not a multiple of $6$. If $a$ is odd we have already seen that $f_a\in {\cal  C}_3$ and so $a\not\in {\cal A}$. If $a$ is even then $a=6\ell \pm 2$ for some $\ell$. So, let $a=6\ell+2$ and $L=a-1$. The previous argument still works the same way because $\gcd(a-2,a+1)=3$ and so we avoid the situation $f_a\in {\cal ID}$. Finally, if $a=6\ell+4$ and $S=a+1$, then $\gcd(a+2,a-1)=3$ and again we avoid $f_a\in {\cal ID}$ as in the case above. In each case  $a\not\in {\cal A}$, showing the other inclusion, i.e. ${\cal A}\subset B$. \end{proof}

For the rest of the paper we are simply using $f$ for $f_3$ and the results and conjecture are going to be concerned with this case.

\begin{theorem}\label{theorem5}%Theorem 2.7

 Let $\p$ denote the $n^{th}$ primorial number. Then for $n\geq 2$, $\p \pm 1$ and $2\p \pm 1 $ are records. 

\end{theorem}

\begin{proof}  If $q=\p -1$ is a prime then we are done because every prime is a record. So, we may assume that  $q$ is composite.

Let $r_l$ denote the largest record less than $q$. Then the next record will be the smallest number larger than $r_l$ but coprime to $r_l-1$. If $q$ is not a record then $r_l-1$ must have a common factor $d$ with $q$. Now $q$ does not have $p_1,p_2,\dots ,p_n$ as prime factors so $d\geq p_{n+1}$. Hence $q-(r_l-1) \geq p_{n+1}$ or $r_i-1+p_{n+}\le q$. Now $r_l-1$ cannot have each of $p_1,p_2,\dots p_n$ as prime factors as it is smaller than $q$. As a result, it must not be divisible by some $p_k$ for $k\leq n$.  Then $r_l-1+p_k$ is a record as it is coprime with $r_l-1$. But $r_l-1+p_k$ lies between $r_l$ and $r_l+p_n<r_l-1+p_{p+1}\le q$ which contradicts the maximality of $r_l$. Hence $q$ must be a record. 

Now, since $f(\p -1)=\p-2$ and $3\mid \p$, $\p -2$ and $\p+1$ must be  coprime. This implies that $\p +1$ is also a record. 

The other part of the theorem's statement follows because there is always a prime between $\p +1$ and $2\p $ which will always be a record. So the largest record less than $2\p -1$, say  $r_{l'}$, will be larger than $\p +1$. Hence $r_{l'}-1$ is larger than $\p$ and less than $2\p$. This means that  $r_{l'}-1$ cannot be divisible by each of the first $n$ primes and a similar argument as above applies. 
\end{proof}

We notice that this shows that $\{\p-1,\p+1\}$  and $\{2\p-1,2\p+1\}$ are twin records. This fact is exploited in the next theorem. 

%%%%%%%%%%%%%%%%%%%%%%%%%%%
\begin{theorem}\label{theorem6}%Theorem 2.8 
\begin{equation}\label{ap}
f(\p+k)=f(k)+\p
\end{equation}
for $k\in [p_{n+1},2\p ]$.
\end{theorem}

\begin{proof} We shall first prove (\ref{ap})  for $k\in [\p ,2\p]$. It suffices to establish a one-to-one correspondence between the records in $I_1:=[\p,2\p ]$ and the ones in $I_2:=[2\p,3\p ]$, i.e.,  show that $r$ is a record in  $I_1$ if and only if $r+\p $ is a record in  $I_2$. This is because between records the values taken by $f$ are determined by the records and their corresponding turning points. The values in between are constructed in a pattern that is compatible with the translation by $\p$. So, $$\p+1=r_1<r_2<r_3<\dots <r_s$$ is all of the records in $I_1$, by Theorem~\ref{theorem5}.

We shall proceed inductively on $r_j$ ($j=1,2,\ldots,s$).  The first record in  $I_1$ is $f(\p)=\p +1$ and the corresponding record in  $I_2$ is $f(2\p)=2\p +1$, by Theorem~\ref{theorem5}.  Now the next record in  $I_1$ will be the smallest number larger than $\p+1$ but coprime to $\p $ which is $r_2=\p +p_{n+1}$ and the next one in $I_2$ will be $2\p+p_{n+1}$.  This is true because there is always a prime $p'$ between $\p $ and $2\p$ which is a record hence $r _2$  has to lie in the interval $I_1$.   (This shows that $p_{n+1} <\p$ which we will use later). We proceed like this to the point where $r_m +\p$ is a record in $I_2$ and $r_m$ is a record in $I_1$ $(m<s)$.  We let $r_m=\p +r$, with $r<\p $.  Then $r_{m+1}$ is the smallest number coprime to $\p +r-1$ and bigger than $\p +r$. Now $r-1$ cannot be divisible by all the primes $p_1,p_2,\dots p_{n}$. Let $p_i$ $(i\leq n)$ be the smallest prime not dividing $r-1$. As $p_i\mid \p $ the smallest number coprime to $\p +r-1$ and bigger than $\p+r$ is $\p+r-1+p_i=r_m-1+p_i$. A similar argument shows that $2\p +r-1+p_i=r_m-1+\p+p_i$ is the next record in $I_2$. Not only that, but every record in $I_2$ appears as a translatsssion of the corresponding record in $I_1$. 

Let us define next $I=[p_{n+1},\p]$. We have just established that the record of $f$ just after $\p+1 $ is $\p+p_{n+1}$. We also know that $p_{n+1}$ is a record.  Now a proof similar to the one given in the previous paragraph holds, by replacing the $I_1$ with $I$ and $I_2$ with $I_1$.
\end{proof}

\n
\begin{rem} As $p_{n+1}< \p$ we have $2\p +p_{n+1}\in [2\p,3\p]$ the largest record less than $3\p-1$ is at least $2\p +p_{n+1}$. Let the largest record be $r_l$.S o $r_l-1>2\p $ and thus cannot be divisible by all of $p_1,p_2,\dots \p_n$ and so an argument similar to that of Theorem 2.7 gives $3\p \pm 1$ are records. Now once we have this we can prove that for $k$ in $[2\p ,3\p]$, $f(k+\p)=f(k)+\p$ (proof similar to Theorem 2.8). So proceeding similarly we get the following theorem whose proof is clear from this remark.\end{rem}

\begin{theorem}\label{theorem7}%Theorem 2.9
 The numbers $r\p \pm 1$ are records for $r\in \{ 1,2,3,\dots ,p_{n+1}-1 \}$.Further $$f(k+\p)=f(k)+\p$$
 for $k\in [p_{n+1},(p_{n+1}-1)\p ]$.
\end{theorem}

For example, one can check that
$$\ f(30+k)=f(k)+30\ \ \text{for all} \  k\in \ [7,181].$$
$$\ f(210+k)=f(k)+210\ \ \text{for all} \  k\in \ [9,2101].$$
Clearly $f$ is not periodic or additive, but the property above suggests some sort of almost periodicity.

\n \begin{rem}  The multiples of the prime $p_n$ appearing as records tend to show a particular pattern because records are sort of translated over large regions by $\p$ giving a pattern to the multiples of that prime appearing as record. This also suggests that getting a good idea of these multiples of primes can give a sieve-like primality test of a number based on the records of this function. 
\end{rem}

 \noindent
In Theorem~\ref{theorem4} we gave a characterization of the set $\mathcal{A}$. Now, from the material developed afterward we can give a better characterization of the set $\mathcal{A}$ which is included in the following theorem.

\begin{theorem}\label{theorem10}%Theorem 2.10
The set $\cal A$ introduced in \ref{Def 3} has the precise description in terms of primes
\begin{equation}\label{setA}
\begin{array}{c}
\mathcal{A}=\{ 2,4 \}\cup \{ a |a=6k, k\in\mathbb{N},a\neq m \cdot \p +6t,\\ \\\mbox{ where }n>3,\ m,t\in\mathbb{N},1\leq t\leq \left\lfloor \frac{p_{n+1}-2}{6}\right\rfloor \}.
\end{array}
\end{equation}
\end{theorem}

\begin{proof} We say that a natural number $a$ ``nice" if $a=6k, k\in\mathbb{N},\ a\neq \p m+6t,\mbox{ where }n>3,m,t\in\mathbb{N},1\leq t\leq \left\lfloor \frac{p_{n+1}-2}{6}\right\rfloor $. We shall show that a natural number $a$ is ``nice" if and only if there is a record $r$ such that $|r-a|\leq 1$ and then we are done by Theorem~\ref{theorem4}.

Let us assume that $a\in \mathbb{N}$ is ``nice" and consider $n$ such that $\p \leq a<p_{n+1}\#$. If $a=m\p$ for $m=1,2,\dots ,p_{n+1}-1$ then we are done because $m\p +1$ is a record. Also, if $p_{n+1}\equiv 1 \pmod{6}$ and $a=m\p +p_{n+1}-1$ then we are also done because $m\p +p_{n+1}$ is a record. Thus,  we may assume $a\in [m\p +p_{n+1},(m+1)\p-6]$ where $m\in\{ 1,2,\dots ,p_{n+1}-1\}$. Hence, there is a record $r$ such that $|r-a|\leq 1$ if and only if there is a record $r'$ with $|r'-a-m\p |\leq 1$. Now, $p_k\#\leq a-m\p< p_{k+1}\#$ where $k<n$. We can repeat this process until we get $x$ such that $a-x<p_4\#$ and after that one can easily check that for any multiple of $6$ less than $p_4\#$ there is a record $r''$, with $r''=x\pm 1$, which is either one larger or less than $x$.

Now let us assume that $a$ is ``not nice". Then $a=m\p +6t$ where $1\leq t\leq \left\lfloor \frac{p_{n+1}-2}{6}\right\rfloor $. Now if $p_{n+1}\equiv -1 \pmod{6}$ then $a\in [m\p+6,m\p +p_{n+1}-5]$ and as there is no record between $m\p +1$ and $m\p+p_{n+1}$ we must have $a\not\in \mathcal{A}$. Similarly,  if $p_{n+1}\equiv 1 \pmod{6}$ then $a\in [m\p+6,m\p+p_{n+1}-7]$ and for the same reason as before $a\not\in\mathcal{A}$.
\end{proof}

\begin{cor}\label{corolary2}%corollary 2.11
The density of the set of ``not nice" numbers, say $\mathcal{A}'$,  is given by the expression

\begin{equation}\label{eqdensityexeptions}
\sum_{ k\geq 4}\left( \left\lfloor \frac{p_{k+1}-2}{6}\right\rfloor - \left\lfloor \frac{p_{k}-2}{6}\right\rfloor\right) \frac{1}{ p_k\#}.
\end{equation}

\end{cor}

\begin{proof} We observe that the set  $$\mathcal{A}':=\{ a | a=\p m+6t \mbox{ where } n>3,\ m,t\in\mathbb{N} ,1\leq t\leq \left\lfloor \frac{p_{n+1}\#-2}{6}\right\rfloor  \}$$
\n is just the union of the sets

$$\mathcal{A}_{n}':=\left\{ a | a=m p_n\# +6t \mbox{ where } m,t\in\mathbb{N} ,1\leq t\leq \left\lfloor \frac{p_{n+1}\#-2}{6}  \right\rfloor \right\}$$

\n for $n=4,5,\ldots $. These sets are disjoint but because $p_{n+1}\#=p_{n+1}\cdot p_n\#$ and $p_{n+1}>p_n$ each element of the set $\mathcal{A}_{n}'$ with $m\ge p_{n+1}$ and $t\leq \left\lfloor \frac{p_{n+1}\#-2}{6}\right\rfloor$ appears in $\mathcal{A}_{n+1}'$ and in all of the subsequent sets $\mathcal{A}_{k}'$ with $k\ge n+1$.  Hence, we can write $\mathcal{A}'$ us union of disjoint sets 

$$\mathcal{B}_{n}:=\left\{ a | a=m p_n\# +6t \mbox{ where } m,t\in\mathbb{N} ,  \left\lfloor \frac{p_{n}\#-2}{6}  \right\rfloor< t\leq \left\lfloor \frac{p_{n+1}\#-2}{6}  \right\rfloor \right\}.$$

As a result, the density of $\mathcal{A}'$ is the sum of the densities of each 
$\mathcal{B}_{n}$. The density of  $\mathcal{B}_{n}$ is  $$\left ( \left\lfloor \frac{p_{n+1}-2}{6}\right\rfloor - \left\lfloor \frac{p_{n}-2}{6}\right\rfloor\right) \frac{1}{ p_n\#} $$

\n because of the periodicity of the elements in $B_n$  modulo $p_n\#$. This gives the formula (\ref{eqdensityexeptions}).
\end{proof}

%%%%%%%%%%%%%%%%%%%%%%%%%%%
Let us denote the set of all records by ${\cal R}$. Written in non-decreasing order gives essentially the sequence A261271. The definition of A261271 is $a_{n+1} = a_n+p-1$, $a_1=1$ and  $p$ is the smallest prime number that is not a factor of $a_n-1$. The equivalence between two concepts is contained in the Theorem~\ref{theorem1}.  

We are interested in the  following limit

\begin{equation}
\label{limitrecords}
\kappa:= \lim_{n\to \infty}\frac{\#\{r|r\in {\cal R}, r\le n\} }{n}
\end{equation}

\n One can see that we can use the property in  (\ref{ap}) to say more about this limit.  We are going to denote by $s_n$ the number of records between $p_n$ and $p_{n+1}$.  More precisely, we  let  
$$s_n=\#\{r|r\in {\cal R}, p_n\le r<p_{n+1}\}.$$
We notice that $s_1=0$, $s_2=1$, $s_3=1$,\ldots, $s_9=2$, $s_{10}=1$,\ldots, $s_{16}=2$, \ldots .So, this sequence is mostly equal to $1$, and when it is greater than $1$, the difference is the number of composite records between the respective consecutive primes. 

\begin{theorem}\label{theorem8}%Theorem 2.10
 The limit (\ref{limitrecords}) exists and we have 
$$\frac{3}{10}-\sum_{k=4}^{\infty}  \frac{p_{k+1}-p_k}{2p_{k}\#}\le  \kappa \le \frac{3}{10}-\sum_{k=4}^{\infty}  \frac{1}{p_{k}\#}.$$
\end{theorem}

\begin{proof}  Let us first show that (\ref{limitrecords}) exists for a subsequence, namely, the limit

\begin{equation}
\label{limitrecords2}
\lim_{n\to \infty}\frac{\#\{r|r\in {\cal R}, r\le \p+1\} }{\p+1}
\end{equation}

\n exists. To simplify notation, we let  $q_n:=\p+1$ and 
$$w_n= \#\{r|r\in {\cal R}, p_{n+1}\le r \le q_n\}. $$
We observe that $w_1=\#\{r|r\in {\cal R},3 \le r \le 3\}=1$,   $w_2=\#\{r|r\in {\cal R},5 \le r \le 7\}=\#\{5,7\}=2$,  

$$w_3=\#\{r|r\in {\cal R},7 \le r \le 31\}=\#\{7, 11, 13, 17, 19, 23, 25, 29, 31\}=9,$$

\n and so on.   We can use Theorem~\ref{theorem6} and the proof of Theorem~\ref{theorem5}, to see, for instance,
that the records in $(32,61]$ are the records in $R_3:=\{7, 11, 13, 17, 19, 23, 25, 29, 31\}$ translated with $30$. Hence their number is the same as $w_3$. Similarly, we can use the same theorems to conclude that the number of records in  $(62,91]$ is also $w_3$. This extends all the way to the interval $(182, 211]$ and so we can say that 
$$w_4=7\cdot 9-1=62,$$

\n as we can see from the following table of the records in the interval $[7,211]$:
$$ \underset{Table\ 1}{\left[
\begin{array}{ccccccccc}
\boxed{ \color{red} 7} & \boxed{ 11} & \boxed{ 13} & \boxed{ 17} & \boxed{ 19} & \boxed{ 23} & \encircled{\color{blue}25} & \boxed{ 29} & \boxed{ 31} \\
 \boxed{ 37 }&  \boxed{ 41} & \boxed{ 43} & \boxed{ 47} &  \encircled{\color{red} 49} & \boxed{ 53} & \encircled{\color{blue}55} & \boxed{ 59} & \boxed{ 61} \\
 \boxed{ 67} & \boxed{ 71} & \boxed{ 73} & \encircled{\color{red} 77} & \boxed{ 79} & \boxed{ 83} & \encircled{\color{blue}85} & \boxed{ 89} & \encircled{\color{red} 91} \\
 \boxed{ 97} & \boxed{101} & \boxed{103} & \boxed{107} & \boxed{109} & \boxed{113} & \encircled{\color{blue}115} & \encircled{\color{red} 119} & \encircled{121} \\
 \boxed{127} & \boxed{131} & \encircled{\color{red} 133} & \boxed{137} & \boxed{139} & \encircled{143} & \encircled{\color{blue}145} & \boxed{149} & \boxed{151} \\
 \boxed{157} & \encircled{\color{red} 161} & \boxed{163} & \boxed{167} & \encircled{169} & \boxed{173} & \encircled{\color{red} 175} & \boxed{179} & \boxed{181} \\
 \encircled{187} & \boxed{191} & \boxed{193} & \boxed{197} & \boxed{199} & \encircled{\color{red} 203} & \encircled{\color{blue}205} & \encircled{209} & \boxed{211} \\
\end{array}
\right]}$$

We can make this argument in general and derive the formula 
\begin{equation}
\label{recurrenceforw}w_{n+1}=w_np_{n+1}-s_{n+1}.
\end{equation}
Hence, after dividing by $p_{n+1}\#$ we get 

$$\frac{w_{n+1}}{p_{n+1}\#}-\frac{w_{n}}{p_{n}\#}=-\frac{s_{n+1}}{p_{n+1}\#}.$$
Summing up these identities we get 
$$0<\frac{w_{n+1}}{p_{n+1}\#}=\frac{w_{3}}{p_{3}\#}-\sum_{k=3}^n \frac{s_{k+1}}{p_{k+1}\#}\le \frac{w_{3}}{p_{3}\#}-\sum_{k=3}^n \frac{1}{p_{k+1}\#}.$$

\n This shows that $\frac {w_{n}}{p_{n}\#}$ is decreasing (hence convergent) and we have an estimate from above for the limit. 
For the other estimate we can use the obvious inequality $s_n\le \frac{1}{2}(p_{n+1}-p_n)$ giving the maximum number of odd integers in the interval $[p_n,p_{n+1})$. 
The limit in  (\ref{limitrecords}) exists (in general) and it is the same as the one for $\frac{w_{n}}{p_{n}\#}$.
 \end{proof} 

We observe that in Table 1, all the records in red and blue will generate composite records in the subsequent intervals, i.e.,  $[p_n\#,p_{n+1}\#]$. These records are multiples of $7$, one on each column and multiples of $5$ all on column $7$, a total of $15$. There are four records that may or may not turn into primes later on, but it is not clear what happens for:  $121$, $169$, $187$ and $209$.  This leaves us with $63-9-7+1-4=44$ primes. 
%%%%%%%%%%%%%%%%%%%%%%%%%%%

 \n

\n \begin{rem} We can use Theorem~\ref{theorem8} to get numeric bounds on $\kappa $. The series of reciprocal of primorial numbers is a convergent series whose value is known to be an irrational number (see \cite{primorial} and \cite{irrat}). The first few terms of the decimal expansion are given by $0.7052301717918009\ldots$. Hence,  we can write $0.704\leq \prod_{n=1}^{\infty}\frac{1}{\p}\leq 0.706$.
This gives us $$\kappa  \leq 0.3-0.704+\frac{1}{2}+\frac{1}{6}+\frac{1}{30}=0.296=\frac{296}{1000}.$$\end{rem}

\noindent Also as $p_{k+1}\leq 2p_k$ , we have $\kappa \geq 0.3-\prod_{n=3}^{\infty}\frac{1}{\p}\geq 0.3-0.706+\frac{1}{2}+\frac{1}{6}=\frac{782}{3000}\approx 0.26067$.

Using the Prime Number Theorem we can tell what is the asymptotic density of primes within the number of records. 

\begin{cor}\label{theorem11}%Theorem 2.11
 We have the following limit
$$ \lim_{n\to \infty}\frac{\# \{r\in {\cal R}|r\ \ \text{is prime and}\ r <n\} }{\# \{r\in {\cal R} |  r<n \} }\ln n=\frac{1}{\kappa} .$$
\end{cor} 

This ratio is included in the next figure for the first $1000$ records. 

\begin{figure}[H]
\centering
\includegraphics[height=0.4\textwidth,width=0.6\textwidth]{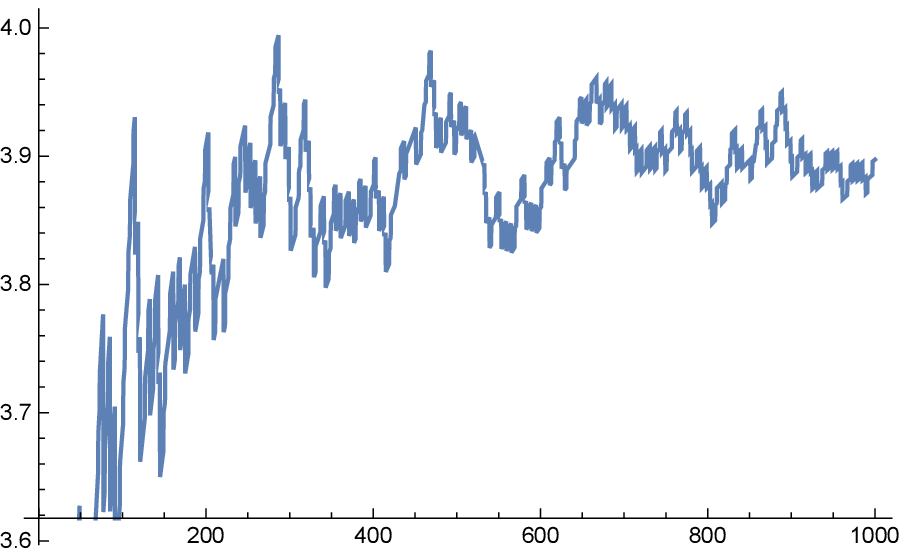}
    \caption{$\frac{\# \{r|r\ \ prime r<n\} }{\# \{r|  r<n \} }\ln n$ }
\end{figure}

\begin{figure}[H]
\centering
\includegraphics[height=0.4\textwidth,width=0.6\textwidth]{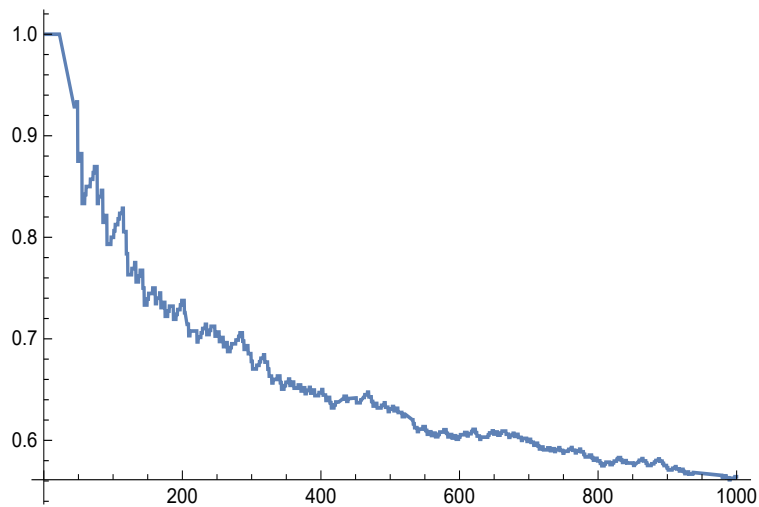}
    \caption{Number of primes within the records}
\end{figure}

\section{Conjectures, questions, problems  and other connections}
\begin{enumerate}[label=\roman*)]

\item We defined in the Introduction the function $h(j)=C(M_{j+1})-C(m_j)$, $j\in \mathbb N$,  where $(m_j,M_j)$ is a twin prime pair (or we can simply refer to only records that form a pair). Numerical evidence suggests that the range of $h$ is $\{k\in\mathbb Z|k\ge 0\}$. So, we are making this conjecture which relates to the twin prime conjecture.

\item Is there a good description for the set  $\overline{{\cal R}}$, the subset of records which are composite ?

We know that all the prime numbers appear as records in $f_3$. So it is important to
characterize the composite numbers that appear as a record in $f_3$, to differentiate the primes.
For that,  let us look at the first few multiples of $5$ that appear as a record:
$$25, 55, 85, 115, 145, 175, 205, 235, 265, 295, 325, 355, 385, 415, 445, 475, \ldots $$
Note that consecutive pairs of these differ by $30$.

Next, let us look at the multiples of $7$ appearing
as a record. They seem to be more interesting:
$$49, 77, 91, 119, 133, 161, 175, 203, 259, 287, 301, 329, 343, 371, 385, 413, 469,\ldots$$
Here the multiples occur on an interval of $28$ and $14$ alternately until
$203$ where it takes a jump of $56$. Then it continues this pattern taking alternate jumps of $28$ and $14$ until in $413=203+210$  where it again takes a jump
of $56$. We guess that this pattern continues until $623=413+210$ where it again
takes a jump of $56$. So these multiples of $7$ appearing as records seem to be very
predictable. Similar observations can be made regarding other primes.  The first
few multiples of $11$ appearing as records are given below:
$$55, 77, 121, 143, 187, 209, 253, 319, 341, 385, 407, 451, 473, 517, 539, 583,\ldots$$
Here the jumps are alternately 22 and 44.
\end{enumerate}
\section{An ad hoc proof of the surjectivity of $f_a$}

Let us show that every prime $p$ is in the range of $f$. If $a=p$ we are done. Otherwise, let $i$ be the greatest index such that $f(i)<p$. This index exists because $f$ is one-to-one (so for $M>0$ there exits $n$ large enough that $f(m)\ge M$  for all $m>n$), and $f(1)=1<p$. If $i=1$ then it must be the case that $1<p<a=f(2)$ and none of the numbers $2$, $3$, \ldots, $p-1$ appear in the sequence $\{f(n)\}_n$ after $a$ by the definition of $i$. Hence, if $\gcd(a,p)=1$ we have $f(3)=p$ and we are done. If not, $a=pa'$ and then $f(3)$ must not contain $p$ in its prime factorization. But then $\gcd(f(3),p)=1$ which forces $f(4)=p$ because of the assumption on $i$. Therefore we may assume $i\ge 2$ and we can apply the definition on $f(i+1)$. Namely, if $p$ is in the set $\{f(1),f(2),\ldots,f(i)\}$ we are done. If not then $f(i+1)=s$ implies $s\ge p$ (by the definition of $i$) and since $\gcd(p,f(i))=1$ we must have $s=p$.

By way of contradiction, we assume that $f$ is not onto. This means there are values in $\mathbb N$ which are not in the range of $f$. Let $k$ be the smallest such number which is not in the range of $f$ (this exists because of the Well-Ordering Principle for  $\mathbb N$).

Then all $m\in \{1,2,\ldots, k-1\}:=A$ must be in the range of $f$. Then the set $B:=f^{-1}(A)=\{f^{-1}(m)|m\in A\}$ has at least $k-1$ elements, but because $f$ is one-to-one  $B$ must have exactly $k-1$ elements. So, we let in order $x_1=1$, $x_2=a$, \ldots, $x_{k-1}=j$ ($x_1<x_2<\ldots<j$). We let  $\ell=f(j)$.  If $\ell = k-1$ then clearly $f(j+1)=k$ and we have a contradiction. A similar argument goes for the situation in which $\ell$ is relatively prime with $k$.

Suppose that the primes in the decomposition of $k$ are $q_1$, $q_2$, $\ldots$, $q_s$ (all distinct primes).

 Then let us look at $\ell_1:=f(j+1)$ which is relatively prime with $\ell$ and since it is not in $A$ we must have $\ell_1>k$. If $\gcd(\ell_1,k)=1$ then by definition we must have $f(j+2)=k$ which is not possible. It remains that $\gcd(\ell_1,k)>1$ so $\ell_1$ and $k$ must have some of the previous primes in common, or in other words at least one of the primes $q_i$ must divide $\ell_1$.

By induction, we can show that $\ell_n:=f(j+n)$ is then a number that must have some prime factor $q_i$ for every $n\ge 1$. This is in contradiction with the fact that all primes must be in the range of $f$.

\section{Acknowledgements}
We thank Mangesh B. Rege who proposed a special case of the problem, as a research question, to the first author while he was in high school. Many thanks are due to Thomas Merino, who wrote a Python code to compute the sequence $f_a(n)$ for significantly big values of $n$ and observed many conjectures about the braking
points and the primorial numbers. Finally, we thank Andrew Ionascu who turned our
sequence into a three-voice canon mdi wave file (\href{https://ejionascu.ro/muzica/andrew2.wav}{canon}).

\end{document}